\documentclass[12pt]{article}

\usepackage{amssymb}
\usepackage{hyperref}
\usepackage{amsmath}

\newtheorem{theo}{\bf \thesection.\arabic{abz}. Theorem }

\newtheorem{prop}{\bf \thesection.\arabic{abz}. Proposition}

\newtheorem{lemm}{\bf \thesection.\arabic{abz}. Lemma}

\newtheorem{coro}{\bf \thesection.\arabic{abz}. Corollary}

\newtheorem{defi}{\bf  \thesection.\arabic{abz}. Definition }

\newtheorem{exa}{\bf \thesection.\arabic{abz}. Example}

\newtheorem{rema}{\bf \thesection.\arabic{abz}. Remark}

\newcounter{abz}[section]
\newcommand{\abz}{\refstepcounter{abz}}

\newcommand\op[1]{\mathop{\rm #1}}

\def\Id{\mathrm{Id}}

\def\F{\mathcal{F}}

\def\R{\mathbb{R}}

\def\P{\mathbb{P}} \def\C{\mathbb{C}}

\def\al{\alpha}

\def\Ga{\Gamma}

\def\d{\partial}

\def\la{\lambda}
\def\La{\Lambda}
\def\om{\omega}

\def\qed{{}\hfill$\square$}

\title{Veronese webs and nonlinear PDEs}

 \author{Boris Kruglikov, Andriy Panasyuk}

\date{}

\begin{document}

\bibliographystyle{plain}

\maketitle

\section{ Introduction}

Three-dimensional Veronese webs are one-parametric foliations of a 3-dimensional space $M$ by surfaces such that their tangents at any point $x$
form a Veronese curve in $\mathop{\rm Gr}_2(T_xM)=\mathbb{P}(T_x^*M)$. They appeared in the study of bi-Hamiltonian systems 
in \cite{GZ}, see also \cite{T3} and the references therein. In \cite{DK} a correspondence between Veronese webs and three-dimensional 
Lorenzian Einstein-Weyl structures of hyper-CR type was established. The latter due to \cite{D} are parametrized by the solutions of the hyper-CR equation
 \begin{equation}\label{hyperCR}
f_{xz}-f_{yy}+f_yf_{xx}-f_xf_{xy}=0.
 \end{equation}
Using the one-to-one correspondence with Veronese webs, the hyper-CR Einstein-Weyl structures were shown by Dunajski and Kry\'{n}ski \cite{DK}
to be also parametrized by the solutions of the dispersionless Hirota equation
 \begin{equation}\label{Hirota}
af_xf_{yz}+bf_yf_{xz}+cf_zf_{xy}=0,\quad a+b+c=0,
 \end{equation}
which was introduced and studied by Zakharevich \cite{Z}.
Both equations above are integrable and the parameters $a,b,c$ are constants, but we will show that they can be taken functions
without destroying the integrability. This can be seen as an integrable deformation, similar to \cite{KM}, though the symmetry is essentially
reduced in this process. The symmetry pseudogroup becomes an equivalence group for the deformation family, which eliminates the functional
parameters though leaves new integrable dispersionless equations.

Here and below by "integrable" we mean those equations that possess a dispersionless Lax pair, and we also show, motivated by \cite{FK},
that they possess Einstein-Weyl structures on solutions, thus representing these equations as reductions of the Einstein-Weyl equation,
integrable by the twistor methods \cite{H}. The introduced equations are not contact equivalent, but they all parametrize Veronese webs and,
using this fact, we will construct B\"{a}cklund transformations between these equations.

The equations that arise are of four types: A, B, C, D. Those that are translationally invariant (the standard requirement for hydrodynamic
integrability test) together with equation \ref{hyperCR} and the universal hierarchy equation \cite{MAS}
are equivalent to the five equations of Ferapontov-Moss \cite{FM} introduced in the context
of quadratic line complexes. Our equations however arise from partially integrable Nijenhuis operators on the way to describe Veronese webs
as a variation of a construction of Zakharevich \cite{Z}.

We establish a correspondence between partially integrable Nijenhuis operators to the operator fields with vanishing Nijenhuis tensor,
and deduce from this a realization of Veronese webs through solutions of equations of any type A, B, C, D. We compute several examples of
realizations, which also provide some exact solutions to the corresponding dispersionless PDEs.

We perceive that the above correspondence can be used as a link between bi-Hamiltonian finite-dimensional integrable systems and
dispersionless integrable PDE related to the Veronese webs. In particular, a combination of results of \cite{PZ} with Theorem \ref{mainTh} suggests
that any generic bi-Hamiltonian system on an odd-dimensional manifold (in particular, necessarily having the Poisson pencil consisting of  degenerate Poisson structures)  can be viewed as a reduction of a bi-symplectic bi-Hamiltonian system of any possible type, i.e.\ whose Poisson pencil is generated by two symplectic forms $\omega_1,\omega_2$ on an even-dimensional manifold such that the Nijenhuis tensor $\omega_2^{-1}\circ\omega_1$ attains all possible types. This application will be addressed in a further research.

\section{Veronese webs and partial Nijenhuis operators}

\abz\label{1}
\begin{defi}\rm {\em A Veronese web} is a collection (below the projective line $\P^1=\R P^1=\R\cup\{\infty\}$, but it can be also taken complex,
i.e., $\C P^1$)
$$
\{\F_\la\}_{\la\in\P^1}
$$
of foliations $\F_\la$ of codimension 1 on $M^{n+1}$ such that
$$
\forall x\in M\ \exists\ \mbox{a\ local\ coframe}\ (\al_0,\ldots,\al_n), \al_i\in\Ga(T^*M)
$$
with
$$
(T_x\F_\la)^\bot=\langle\al_0+\la\al_1+\cdots+\la^n\al_n\rangle.
$$
\end{defi}

In particular, $(T_x\F_0)^\bot=\langle\al_0\rangle$ and $(T_x\F_\infty)^\bot=\langle\al_n\rangle$. Moreover, the frame $\al$ is defined up to
simultaneous multiplication by a nonvanishing smooth function.

\abz\label{2}
\begin{defi}\rm
A {\it partial Nijenhuis operator} ({\em PNO}) on a manifold $M$ is a pair $(\F,\bar{J})$, where $\F$ is a foliation on $M$ and
$\bar{J}:T\F\to TM$ is a partial (1,1)-tensor such that $\forall\ X,Y\in \Ga(T\F)$
\begin{enumerate}\item $[X,Y]_{\bar{J}}:=[\bar{J}X,Y]+[X,\bar{J}Y]-\bar{J}[X,Y]\ \in \ \Ga(T\F)$;
\item $N_{\bar{J}}(X,Y):=[\bar{J}X,\bar{J}Y]-\bar{J}[X,Y]_{\bar{J}}=0$.\end{enumerate}
\end{defi}

\abz\label{3}
\begin{exa}
\rm Let $J:TM\to TM$ be a {\it Nijenhuis operator}, i.e., a (1,1)-tensor on $M$ such that its Nijenhuis tensor vanishes, $N_J\equiv 0$.
Then $(M,J)$ is a PNO.
\end{exa}

The following statements are straightforward.

\abz\label{4}
 \begin{lemm}
Let $(\F,\bar{J})$ be a PNO on $M$. Then
 \begin{itemize}
\item $(\F,\bar{J}_\la)$ is a PNO; here $\bar{J}_\la:=\bar{J}-\la I$, and $I:T\F\hookrightarrow TM$ is the canonical inclusion;
\item $[X,Y]_{\bar{J}_\la}$ is a Lie bracket on $\Ga(T\F)$;
\item $\bar{J}_\la:\Ga(T\F)\to \Ga(TM)$ is a homomorphism of Lie algebras.
 \end{itemize}
In particular, if $\bar{J}_\la(T\F) \subset TM$ is a distribution, it is integrable:
 $$
\hspace{7.4cm}\bar{J}_\la(T\F)=T\F_\la. \hspace{7.4cm}\square
 $$
 \end{lemm}

\abz\label{5}
 \begin{lemm}
Let $J:TM\to TM$ be a Nijenhuis operator ($N_J\equiv 0$), $\F$ a foliation, and assume
 \begin{itemize}
\item $\bar{J}:=J|_{T\F}:T\F\to TM$ is injective;
\item $J(T\F)\subset TM$ is an integrable distribution.
 \end{itemize}
Then $(\F,\bar{J})$ is a PNO. \qed
\end{lemm}

\bigskip

\abz\label{rem0}
\begin{rema}\rm
The converse to Lemma \ref{5} will be obtained in Theorem \ref{mainTh}.
Notice though that a Nijenhuis operators $J:TM\to TM$ inducing the given PNO $(\F,\bar{J})$ can be non-unique.
\end{rema}

\abz\label{rem01}
\begin{rema}\rm The notion of partial Nijenhuis operator is elaborated in detail in \cite{PZ} (in particular Lemmas \ref{4}, \ref{5} are borrowed from there). Similar notions appeared in \cite{cgm} (under the name ``outer Nijenhuis tensor'') and \cite{T4} (without name and in dual terms, i.e. using differential forms).

\end{rema}

Let us say that a PNO $(\F,\bar{J})$ is of \emph{of generic type} if the pair of operators $\bar{J},I:T\F\hookrightarrow TM$ has unique Kronecker block
in their Jordan--Kronecker decomposition, i.e., there exist local frames $v_1,\ldots,v_n\in\Ga(T\F)$, $w_0,\ldots,w_n\in\Ga(TM)$, in which
$$
\bar{J}=\left[
    \begin{array}{cccc}
      1 & & &  \\
      0 & 1 &  &\\
      & \ddots&\ddots& \\
       & & 0 & 1 \\
       & & & 0 \\
    \end{array}
  \right],
I=\left[
    \begin{array}{cccc}
      0 & & &  \\
      1 & 0 &  &\\
            & \ddots&\ddots& \\
       & & 1 & 0 \\
       & & & 1 \\
    \end{array}
  \right].
$$

\abz\label{theo1}
\begin{theo} There exists a 1-1-correspondence between Veronese webs $\{\F_\la\}$ on $M^{n+1}$ and PNOs $(\F,\bar{J})$ of generic type
such that $\F_\infty=\F$ and $T\F_0=\bar{J}\,T\F$. Locally one can choose $\bar{J}=J|_{T\F}$ for a Nijenhuis operator $J$ as in Lemma \ref{5}.
\end{theo}

\noindent {\sc Proof.} ($\Longleftarrow$) $(\F,\bar{J})\mapsto \bar{J}_\la(T\F)=T\F_\la$ (use Lemma \ref{4}).

($\Longrightarrow$)\emph{Variation of a construction of F.J.\ Turiel \cite{T1}:}

Let $\{\F_\la\}$ be a Veronese web on $M^{n+1}$. Fix pairwise distinct nonzero numbers $\la_1,\dots,\la_{n+1}$. Then
$$
D_i=\cap_{j\neq i}T\F_{\la_j},\quad i=1,\dots,n+1,
$$
are 1-dimensional distributions such that $D_i+D_j$ are integrable rank 2 distributions for all $i\neq j$.
Hence there exists a local coordinate system $(x_1,\dots,x_{n+1})$ such that $D_i=\langle \partial_{x_i}\rangle$. Define
$$
J\partial_{x_i}=\la_i\partial_{x_i}.
$$
Then $N_J\equiv 0$ and $(\F_\infty,\bar{J})$ with $\bar{J}=J|_{T\F_\infty})$ is a PNO. Indeed $J(T\F_\infty)=T\F_0$ is integrable and one can use
Lemma \ref{5}. Moreover, $J_{\la_i}(T\F_\infty)=T\F_{\la_i}$, $i=1,\ldots,n+1$ and by the uniqueness property of the Veronese curve $J_\la(T\F_\infty)=\F_\la$.

The constructed PNO $(\F_\infty,\bar{J})$ is independent of the choice of the numbers $\la_i$. Indeed, let
$(T\F_\la)^\bot=\langle \al_0+\la \al_1+\dots+\la^n\al_n\rangle=:\langle\al^\la\rangle$ and let $X_0,\dots,X_n$ be the frame
dual to the coframe $\al_0,\dots,\al_n$. Then the partial operator $\bar{J}: T\F_\infty=\langle X_0,\dots,X_{n-1}\rangle\to TM$
satisfying $\al^\la(\bar{J}_\la\,T\F_\infty)=0$ for any $\la$ is uniquely determined by $\bar{J}X_k=X_{k+1}$, $0\le k<n$.
Note also that the pair $(\bar{J},I)$ has canonical Jordan--Kronecker matrix form in the frames $X_0,\dots,X_{n-1}$ and $X_0,\dots,X_n$.
\qed

\section{The Hirota equation}
\label{Hiro}

Starting from this section we assume $\dim M=3$. Some results hold for the general dimension $n$, but for simplicity we assume $n=2$.
We begin with the following
\smallskip

\emph{Variation of a construction of I.\ Zakharevich \cite{Z}:}\\
It follows from Theorem \ref{theo1} that, given a Veronese web, one can construct a PNO that, at least locally, can be extended to a
Nijenhuis operator defined on the whole tangent bundle $TM$. Conversely, starting from a (1,1)-tensor $J$ we want to construct a PNO $(\F,\bar{J})$,
$\bar{J}=J|_{T\F}$, which corresponds to a Veronese web by Theorem \ref{theo1}. Assuming that the foliation $\F$ is given by
$f=\op{const}$ for some smooth function $f$ and that $J$ is semi-simple, we will use Lemma \ref{5} to obtain sufficient conditions for $\bar{J}$
to be a PNO in terms of a PDE on $f$ in which we recognize the Hirota equation (\ref{Hirota}).

Consider $M=\R^3(x_1,x_2,x_3)$ and let $\la_1,\la_2,\la_3$ be pairwise distinct nonzero numbers.
Construct a Nijenhuis operator $J:T\R^3\to T\R^3$ by
 \begin{equation}\label{Jx}
J\partial_{x_i}=\la_i\partial_{x_i}.
 \end{equation}
Let $f:\R^3\to \R$ be a smooth function such that $f_{x_i}\not=0$. Define a foliation $\F_\infty$ by $f=\op{const}$, i.e.,
$T\F_\infty:=\langle df\rangle^\bot$. Then $(J(T\F_\infty))^\bot=\langle\om\rangle$, where
$$
\om= \la_1^{-1}f_{x_1}dx_1+\la_2^{-1}f_{x_2}dx_2+\la_3^{-1}f_{x_3}dx_3.
$$
The distribution $J(T\F_\infty)$ is integrable and if and only if  $d\om\wedge\om=0$, i.e.,
 \begin{equation}\label{Hirota-a}
(\la_2-\la_3)f_{x_1}f_{x_2x_3}+(\la_3-\la_1)f_{x_2}f_{x_3x_1}+(\la_1-\la_2)f_{x_3}f_{x_1x_2}=0.
 \end{equation}
The following theorem is a variant of \cite[Theorem 3.8]{Z} (our proof is different).

\abz\label{6}
\begin{theo} Let $\la_1,\la_2,\la_3$ be  distinct real numbers.
 \begin{enumerate}
\item For any solution $f$ of (\ref{Hirota-a}) on a domain $U\subset M$ with $f_{x_i}\not=0$ for all $i=1,2,3$ the 1-form
 \begin{equation}\label{1-f}
\al^\la=(\la_2-\la)(\la_3-\la)f_{x_1}dx_1+(\la_3-\la)(\la_1-\la)f_{x_2}dx_2+(\la_1-\la)(\la_2-\la)f_{x_3}dx_3
 \end{equation}
defines a Veronese web $\F_\la$ on $U$ by $T\F_\la=\langle\al^\la\rangle^\bot$. We have:
 \begin{equation}\label{FFF}
\F_{\la_i}=\{x_i=const\},\ \F_\infty=\{f=const\}.
 \end{equation}
\item Conversely, let $\F_\la$ be a Veronese web on a 3-dimensional smooth manifold $M$. Then in a neighbourhood of any point on $M$ there exist local coordinates $(x_1,x_2,x_3)$ such that any smooth first integral $f$ of the foliation $\F_\infty$ is a solution of equation (\ref{Hirota-a}) with $f_{x_i}\not=0$.
 \end{enumerate}
Consequently, we obtain a 1-1-correspondence between Veronese webs $\F_\la$ satisfying (\ref{FFF}) and
the classes $[f]$ of solutions $f$ of (\ref{Hirota-a}) with  $f_{x_i}\not=0$ modulo the following equivalence relation:
$f\sim g$ if there exist local diffeomorphisms $\psi,\phi_1,\phi_2,\phi_3$ of $\R$ that $f(x_1,x_2,x_3)=\psi(g(\phi_1(x_1),\phi_2(x_2),\phi_3(x_3))$
(obviously, if $f\sim g$ and $f$ solves (\ref{Hirota-a}), then $g$ does the same).
 \end{theo}

\noindent {\sc Proof.}
On a solution $f$ of equation (\ref{Hirota-a}) we get $d\om\wedge\om=0$, and so the distribution $J(T\F_\infty)$ is integrable.
Consequently, $\bar{J}=J|_{T\F_\infty}$ is a PNO by Lemma \ref{2}. The condition $f_{x_i}\not=0$
implies that the pair $(\bar{J},I)$ has generic type and thus defines a Veronese web $\F_\la$ by Theorem \ref{theo1}.
The Veronese curve $\al^\la$ in $T^*U$ such that $(T\F_\la)^\bot=\langle \al^\la \rangle$ annihilates the distribution $\bar{J}_\la(T\F_\infty)=T\F_\la$.
Direct check shows that it is given by formula (\ref{1-f}), in particular satisfies (\ref{FFF}).

Conversely, let $\F_\la$ be a Veronese web and $f$ a first integral of $\F_\infty$. The proof of Theorem \ref{theo1} yields the coordinates
$(x_1,x_2,x_3)$ such that (\ref{Jx}) holds. The distribution $J(T\F_\infty)=T\F_0$ is integrable, hence $d\om\wedge\om=0$ and $f$
solves (\ref{Hirota-a}). The condition $f_{x_i}\not=0$ follows from nondegeneracy of $\al^\la$.

Finally, the last statement follows from the fact that the first integrals of the three Veronese foliations corresponding to different $\la_1,\la_2,\la_3$
determine the first integral of any other foliation up to postcomposition with a local diffeomorphism.
 \qed

\bigskip

The following vector fields depending on the parameter $\la$ annihilate the 1-form $\al^\la$ and is a Lax pair for the Hirota equation (\ref{Hirota-a}):
 $$
v^\la:=(\la-\la_1)f_{x_2}\partial_{x_1}-(\la-\la_2)f_{x_1}\partial_{x_2}, \quad
w^\la:=(\la-\la_2)f_{x_3}\partial_{x_2}-(\la-\la_3)f_{x_2}\partial_{x_3}.
 $$

\section{Classification of Nijenhuis operators in 3D}\label{3DNO}

We want to extend the construction of the previous section using other PNO $(\F,\bar{J})$, where
$\bar{J}=J|_{T\F}$ is the restriction of a Nijenhuis operator. For this we need to describe the Nijenhuis operators in 3D.

Let us call a germ of a (1,1)-tensor (operator field) {\em stable} if its Jordan normal form does not bifurcate at this point, and the multiplicities
of the eigenvalues do not change in a neighborhood. Let us call it {\em non-degenerate} if no eigenvalues corresponding to different Jordan blocks
are equal. While the first assumption reduces complications with classification of Nijenhuis operators (only a finite typical germ of non-stable Nijenhuis
operator can be classified via singularity theory), the second assumption removes degenerate PDE (that do not produce Veronose webs),
corresponding to a Nijenhuis operator, so we adapt both assumptions.

\abz\label{theoB1}
\begin{theo}
A germ of stable non-degenerate Nijenhuis operator in $\R^3$ has one of the four possible forms in a local coordinate system $(x_1,x_2,x_3)$:
 \begin{itemize}
\item[{\rm({\bf A})}] Real semi-simple case:
$J=\left[\begin{array}{ccc} \la_1(x_1) & 0 & 0\\ 0 & \la_2(x_2) & 0\\ 0 & 0 & \la_3(x_3) \end{array}\right]$,
\item[{\rm({\bf B})}] $2\times 2$ and $1\times1$ Jordan blocks:
$J=\left[\begin{array}{ccc} \la_2(x_2) & 1 & 0\\ 0 & \la_2(x_2) & 0\\ 0 & 0 & \la_3(x_3) \end{array}\right]$,
\item[{\rm({\bf C})}] $3\times 3$ Jordan block:
$J=\left[\begin{array}{ccc} \la_3(x_3) & e^{\la_3'(x_3)x_2} & 0\\ 0 & \la_3(x_3) & 1\\ 0 & 0 & \la_3(x_3) \end{array}\right]$,
\item[{\rm({\bf D})}] Complex semi-simple case:
$J=\left[\begin{array}{ccc} a(x_1,x_2) & -b(x_1,x_2) & 0\\ b(x_1,x_2) & a(x_1,x_2) & 0\\ 0 & 0 & \la_3(x_3) \end{array}\right]$.
 \end{itemize}
Here $\la_i(x_i)$ are arbitrary smooth functions ($\la_i\neq\la_j$ if $i\neq j$), and $Z(z)=a(x_1,x_2)+i\,b(x_1,x_2)$, $z=x_1+i\,x_2$, is an arbitrary holomorphic function
($a$ and $b$ are harmonic duals).
\end{theo}

\noindent {\sc Proof.}
Notice that if $J v_i=\la_i v_i$, $J v_j=\la_j v_j$ for vector fields $v_i,v_j$ ($\la_i$ are functions), then
 $$
N_J(v_i,v_j)=(J-\la_i)(J-\la_j)[v_i,v_j]+(\la_i-\la_j)(v_i(\la_j)v_j+v_j(\la_i)v_i)=0.
 $$
When the spectrum $\op{Sp}(J)$ is real simple, this implies that the distributions $\langle v_i,v_j\rangle$ are integrable and $v_i(\la_j)=0$ for $i\neq j$.
Hence the eigendistributions are jointly integrable, giving the coordinate system $(x_1,x_2,x_3)$ in which
 \[
J=\la_1(x_1)\,\partial_{x_1}\otimes dx_1 + \la_2(x_2)\,\partial_{x_2}\otimes dx_2 + \la_3(x_3)\,\partial_{x_3}\otimes dx_3. \tag{A}\label{J-A}
 \]

In the case $\op{Sp}(J)$ is complex simple, the condition $N_J=0$ implies integrability of eigendistributions, and so our space is locally a product
of $\C^1(z)=\R^2(x_1,x_2)$ and $\R(x_3)$. Now one easily checks that $N_J=0$ implies the Cauchy-Riemann equations on $Z=a+ib$, and we conclude
 \[
J=a(x_1,x_2)\,(\partial_{x_1}\otimes dx_1+\partial_{x_2}\otimes dx_2) + b(x_1,x_2)\,(\partial_{x_2}\otimes dx_1-\partial_{x_1}\otimes dx_2) +
 \la_3(x_3)\,\partial_{x_3}\otimes dx_3. \tag{D}\label{J-D}
 \]

Consider now the case of $2\times 2$ Jordan block with eigenvalue $\la_2$ and $1\times1$ block with eigenvalue $\la_3$ (recall $\la_i$ are functions).
Thus there exists a frame $(v_1,v_2,v_3)$ such that
 $$
J v_1= \la_2 v_1,\quad J v_2=\la_2 v_2+v_1,\quad J v_3=\la_3 v_3.
 $$
We compute:
 \begin{gather*}
N_J(v_1,v_2)=(J-\la_2)^2[v_1,v_2]-2v_1(\la_2)v_1=0;\\
N_J(v_1,v_3)=(J-\la_2)(J-\la_3)[v_1,v_3]+(\la_2-\la_3)(v_1(\la_3)v_3+v_3(\la_2)v_1)=0;\\
N_J(v_2,v_3)=(J-\la_2)(J-\la_3)[v_2,v_3]-(J-\la_3)[v_1,v_3]+(\la_2-\la_3)(v_2(\la_3)v_3\\ \hspace{7cm}+v_3(\la_2)v_2)+v_1(\la_3)v_3+v_3(\la_2)v_1=0.
 \end{gather*}
The first equation implies that the distribution $\langle v_1,v_2\rangle$ is integrable, and $v_1(\la_2)=0$. The second yields $v_1(\la_3)=0$,
and the third $v_2(\la_3)=0$.
Applying $(J-\la_2)$ to the last equation and comparing the result to the second we get integrability of the distribution $\langle v_1,v_3\rangle$
and $v_3(\la_2)=0$.

Notice that we can freely change $v_2$ by $v_1$, so we can arrange that all three distributions $\langle v_i,v_j\rangle$ are integrable, whence we
get the coordinate system $(x_1,x_2,x_3)$ in which
 $$
J\partial_{x_1}=\la_2(x_2)\partial_{x_1},\quad J\partial_{x_2}=\la_2(x_2)\partial_{x_2}+\nu\partial_{x_1},\quad
J\partial_{x_3}=\la_3(x_3)\partial_{x_3}
 $$
and $\nu\neq0$ is some function that, due to the condition $N_J=0$, satisfies $\nu_{x_3}=0$.

Thus $J$ restricts to $\R^2(x_1,x_2)=\{x_3=\op{const}\}$ planes, and these restrictions do not depend on $x_3$. There exists
(by the Hamilton-Jacobi theory of first order PDE) a function $\kappa$ solving the constraint $[\nu\partial_{x_1},\partial_{x_2}+\kappa\partial_{x_1}]=0$.
Then we can take new coordinates in $\R^2(x_1,x_2)$ such that the above pair of vector fields is $\partial_{x_1},\partial_{x_2}$.
In these coordinates $\nu=1$, and we obtain
 \[
J=\la_2(x_2)\,(\partial_{x_1}\otimes dx_1+\partial_{x_2}\otimes dx_2) + \partial_{x_1}\otimes dx_2+
\la_3(x_3)\,\partial_{x_3}\otimes dx_3. \tag{B}\label{J-B}
 \]

Finally, consider the Jordan block of size $3\times 3$ with eigenvalue $\la_3$, i.e., there exists a frame $v_1,v_2,v_3$ such that
 $$
J v_1=\la_3 v_1,\quad J v_2= \la_3 v_2+ v_1,\quad J v_3= \la_3 v_3+ v_2.
 $$
We compute:
 \begin{gather*}
N_J(v_1,v_2)=(J-\la_3)^2[v_1,v_2]-2v_1(\la_3)v_1=0;\\
N_J(v_1,v_3)=(J-\la_3)^2[v_1,v_3]-(J-\la_3)[v_1,v_2]-(v_1(\la_3)v_2+v_2(\la_3)v_1)=0;\\
N_J(v_2,v_3)=(J-\la_3)^2[v_2,v_3]-(J-\la_3)[v_1,v_3]+[v_1,v_2]+v_1(\la_3)v_3+v_3(\la_3)v_1-2v_2(\la_3)v_2=0.
 \end{gather*}
Applying $(J-\la_3)$ to the second equation and comparing to the first we conclude that the distribution $\langle v_1,v_2\rangle$ is integrable
and $v_1(\la_3)=0$. Applying $(J-\la_3)$ to the third equation and comparing to the second we conclude that $v_2(\la_3)=0$.
The flag $\langle v_1\rangle\subset \langle v_1,v_2\rangle\subset TM$ is rectifyable.

Thus we can find coordinates $(x_1,x_2,x_3)$ such that for some functions $\nu\neq0$, $\sigma\neq0$ and $\tau$
 $$
J\partial_{x_1}=\la_3(x_3)\partial_{x_1},\quad J\partial_{x_2}=\la_3(x_3)\partial_{x_2}+\nu\partial_{x_1},\quad
J\partial_{x_3}=\la_3(x_3)\partial_{x_3}+\sigma\partial_{x_2}+\tau\partial_{x_1}.
 $$
The last term with $\tau$ is unavoidable in a general rectification, but since our freedom in the coordinate change is triangular
$(x,y,z)\mapsto(X(x,y,z),Y(y,z),Z(z))$, we can adjust the coordinates to fix $\tau=0$. This is the first normalization.
The condition $N_J=0$ is equivalent to $\sigma_{x_1}=0$, $\nu\,\la'_3(x_3)=(\sigma\nu)_{x_2}$.

The second normalization is to fix $\sigma$. For this notice that by variation $\partial_{x_3}\mapsto\eta=\partial_{x_3}+\kappa\partial_{x_2}$
we can achieve $[\eta,\sigma\partial_{x_2}]=0$ and since $\sigma_{x_1}=0$, we can also choose $\kappa_{x_1}=0$. Then
the basis $(\partial_{x_1},\sigma\partial_{x_2},\eta)$ is holonomic, and changing coordinates to make it the coordinate basis we get $\sigma=1$.

The constraint on $\nu$ now becomes $\nu\,\la'_3(x_3)=\nu_{x_2}$, so that $\nu=b(x_1,x_3)e^{\la'_3(x_3)x_2}$. The third normalization is to fix $b$
(which is nonzero since $\nu$ is nonzero). This is done by the change of coordinates $(x_1,x_2,x_3)\mapsto(\int\frac{dx_1}{b(x_1,x_3)},x_2,x_3)$,
and we obtain
 \[
J=\la_3(x_3)\,(\partial_{x_1}\otimes dx_1+\partial_{x_2}\otimes dx_2+\partial_{x_3}\otimes dx_3) +
e^{\la'_3(x_3)x_2}\partial_{x_1}\otimes dx_2 + \partial_{x_2}\otimes dx_3. \tag{C}\label{J-C}
 \]
This finishes the proof. \qed

\smallskip

\abz\label{remB1}
\begin{rema}\rm
Let us fix the freedom in coordinates for every case of the theorem. These are the corresponding equivalence groups, each being
parametrized by 3 functions of 1 argument.
 \begin{enumerate}
\item[{\rm(A)}] $x_1\mapsto X_1(x_1)$,  $x_2\mapsto X_2(x_2)$,  $x_3\mapsto X_3(x_3)$,
\item[{\rm(B)}] $x_1\mapsto X_2'(x_2)x_1+X_1(x_2)$, $x_2\mapsto X_2(x_2)$, $x_3\mapsto X_3(x_3)$,
\item[{\rm(C)}] $x_1\mapsto X_3'(x_3)\exp\Bigl(\la_3'(x_3)\frac{X_2(x_3)}{X_3'(x_3)}\Bigr)x_1+
\Bigl(\frac{X_3''(x_3)}{\la_3'(x_3)}x_2+\frac{X_2'(x_3)}{\la_3'(x_3)}-\frac{X_3''(x_3)}{\la_3'(x_3)^2}\Bigr)
\exp\Bigl(\la_3'(x_3)\bigl(x_2+\frac{X_2(x_3)}{X_3'(x_3)}\bigr)\Bigr)+X_1(x_3)$,
$x_2\mapsto X_3'(x_3)x_2+X_2(x_3)$, $x_3\mapsto X_3(x_3)$,
\item[{\rm(D)}] $(x_1,x_2)\mapsto X_{12}(x_1,x_2)$,  $x_3\mapsto X_3(x_3)$, where $X_{12}$ is a harmonic map of $\R^2(x_1,x_2)=\C(z)$.
 \end{enumerate}

Every equivalence group acts on the corresponding space of functional parameters, which consists of 3, 2, 1 and 3 functions respectively in cases
A, B, C, D. It is easy to see that every function $\la_i$ from the general stratum of the functional parameter space in the theorem can be
reduced by this equivalence group to either constant or linear function in the corresponding argument. Thus we get the following
normal forms of the Nijenhuis operators near a generic point:
 \begin{enumerate}
\item[{\rm(A-C)}] $\la_i(x_i)=\la_i=\op{const} $ or $\la_i(x_i)=x_i$ for $i=1,2,3$,
\item[{\rm(D)}] $Z(z)=X_{12}(x_1,x_2)$ is either a $\op{const}$ or it equals to $z=x_1+ix_2$; for $\la_3(x_3)$ the same as above.
 \end{enumerate}
Such normal forms are known from the work of Turiel \cite{T}, and we summarize these results in the Appendix (we note there that an additional
assumption of cyclicity taken in loc.cit.\ in the general case is not needed in our 3D case). It is straightforward to see
that all his forms are specifications of ours as just indicated. For instance, case $A_0$ corresponds to $\la_i(x_i)=x_i$, $i=1,2,3$.

The only confusion can come from form $C_0$, which does not resemble Jordan normal form apparent in our case C. Yet that $C_0$ corresponds
to our ${\rm C}|_{\la_3(x_3)=x_3}$. To see this, denote our Nijenhuis operator with specified functional parameters by $J^C_\dagger$
and that of $C_0$ by $J^C_0$, i.e., we have:
 \[
J^C_\dagger=x_3\cdot\op{Id}\nolimits_3+ e^{x_2}\partial_{x_1}\otimes dx_2 + \partial_{x_2}\otimes dx_3,\quad
J^C_0=\xi_3\cdot\op{Id}\nolimits_3+ \partial_{\xi_2}\otimes d\xi_1 + \partial_{\xi_1}\otimes d\xi_3 - \xi_2\,\partial_{\xi_2}\otimes d\xi_3,
 \]
where $\op{Id}_3=\sum_1^3\partial_{x_i}\otimes dx_i=\sum_1^3\partial_{\xi_i}\otimes d\xi_i$. Then for the diffeomorphism
$\Phi:\R^3\to\R^3$, $\Phi(x_1,x_2,x_3)=(\xi_1,\xi_2,\xi_3)$, given by $\xi_1=x_2, \xi_2=x_1e^{-x_2}, \xi_3=x_3$, we get
$\Phi_*J^C_\dagger=J^C_0$.
\end{rema}

\section{PNO deformation of the Hirota and three other PDE}
\label{listPDE}

Using the classification of Nijenhuis operators $J$ from the previous section, we obtain functional families of PNO $(\F,\bar{J})$, $\bar{J}=J|_{T\F}$.
This gives a Veronese web $\F_\la$ with $\F=\F_\infty$ as before: The Veronese curve $\langle\al^\la\rangle$ in $\mathbb{P}T^*M$
is given by the formula $\langle \al^\la \rangle=(T\F_\la)^\bot$, $\bar{J}_\la(T\F_\infty)=T\F_\la$.

Choosing $\om=\op{const}\cdot\al^0$ via $(J(T\F_\infty))^\bot=\langle\om\rangle$,
the Frobenius integrability condition $d\om\wedge\om=0$ written via a first integral $f$ of $\F$ (i.e., $\omega=J^{-1}df$)
is a second order PDE on $f$ involving the functional parameters from $J$. This will be treated as an integrable deformation, cf.\ \cite{KM}.
In loc.cit.\ the deformation was governed by the symmetry algebra, but in our case it is governed by the geometry of  the Nijenhuis tensor.

Just by construction the Veronese web satisfies the condition $\al^\la\wedge d\al^\la=0$, so the obtained equations are integrable
via the Lax pair with spectral parameter $\la$, viz.\ given by the vector fields $v^\la,w^\la$ spanning $(\al^\la)^\perp$.

Let us list the PDEs on the function $f$, corresponding to the cases A, B, C, D of $J$, and indicate the Veronese curves $\al^\la$
(the formulae for $\omega$ and $v^\la,w^\la$ follow).

 \begin{enumerate}
\item[\bf(A)]\hspace{0.5cm}
\framebox{ $(\la_2(x_2)-\la_3(x_3))f_{x_1}f_{x_2x_3}+(\la_3(x_3)-\la_1(x_1))f_{x_2}f_{x_3x_1}+(\la_1(x_1)-\la_2(x_2))f_{x_3}f_{x_1x_2}=0$}
 $$
\hspace{-0.7cm}\al^\la = (\la_2(x_2)-\la)(\la_3(x_3)-\la)f_{x_1}dx_1+
(\la_3(x_3)-\la)(\la_1(x_1)-\la)f_{x_2}dx_2+(\la_1(x_1)-\la)(\la_2(x_2)-\la)f_{x_3}dx_3.
 $$
\item[\bf(B)]\hspace{0.5cm}
\framebox{ $f_{x_1}f_{x_1x_3}-f_{x_3}f_{x_1x_1}+(\la_2(x_2)-\la_3(x_3))(f_{x_1}f_{x_2x_3}-f_{x_2}f_{x_1x_3})
+\la_2'(x_2)f_{x_1}f_{x_3}=0$}
 $$
\hspace{-0.7cm}\al^\la = (\la_2(x_2)-\la)(\la_3(x_3)-\la)(f_{x_1}dx_1+f_{x_2}dx_2)+(\la_2(x_2)-\la)^2f_{x_3}dx_3
-(\la_3(x_3)-\la)f_{x_1}dx_2.
 $$
\item[\bf(C)]\hspace{0.5cm}
\framebox{ $f_{x_1}f_{x_1x_3}-f_{x_3}f_{x_1x_1}+e^{-\la_3'(x_3)x_2}(f_{x_2}f_{x_1x_2}-f_{x_1}f_{x_2x_2})+\la_3''(x_3)x_2f_{x_1}^2=0$}
 $$
\hspace{-0.7cm}\al^\la = (\la_3(x_3)-\la)^2df-(\la_3(x_3)-\la)(e^{\la_3'(x_3)x_2}f_{x_1}dx_2+f_{x_2}dx_3)+e^{\la_3'(x_3)x_2}f_{x_1}dx_3.
 $$
\item[\bf(D)]\hspace{0.5cm}
\framebox{\parbox{15.3cm}{ $a(x_1,x_2)(f_{x_1}f_{x_2x_3}-f_{x_2}f_{x_1x_3})+b(x_1,x_2)(f_{x_3}f_{x_1x_1}+f_{x_3}f_{x_2x_2}-f_{x_1}f_{x_1x_3}-f_{x_2}f_{x_2x_3})\\
\hphantom{A}\hspace{9cm}+\la_3(x_3)(f_{x_2}f_{x_1x_3}-f_{x_1}f_{x_2x_3})=0$}}
 $$
\hspace{-0.7cm}\al^\la = (\la_3(x_3)-\la)((a-\la)(f_{x_1}dx_1+f_{x_2}dx_2)+b(f_{x_1}dx_2-f_{x_2}dx_1))+((a-\la)^2+b^2)f_{x_3}dx_3.
 $$
 \end{enumerate}
Recall that $a=a(x_1,x_2)$ and $b=b(x_1,x_2)$ are harmonic dual, i.e., $a_{x_1}=b_{x_2}$, $a_{x_2}=-b_{x_1}$.

When we come to the normal forms of the Nijenhuis operators at generic points as in Appendix, then we follow Remark \ref{remB1} and obtain the
following specifications of the above equations. Below $c_1,c_2,c_3$ are arbitrary different constants.

 \begin{enumerate}
\item[\bf(A)]
$A_0$: $\la_1(x_1)=x_1, \la_2(x_2)=x_2, \la_3(x_3)=x_3$;\qquad
$A_1$: $\la_1(x_1)=x_1, \la_2(x_2)=x_2, \la_3(x_3)=c_3$;\\
$A_2$: $\la_1(x_1)=x_1, \la_2(x_2)=c_2, \la_3(x_3)=c_3$;\qquad\
$A_3$: $\la_1(x_1)=c_1, \la_2(x_2)=c_2, \la_3(x_3)=c_3$.
\item[\bf(B)]
$B_0$: $\la_2(x_2)=x_2, \la_3(x_3)=x_3$;\qquad
$B_1$: $\la_2(x_2)=x_2, \la_3(x_3)=c_3$;\\
$B_2$: $\la_2(x_2)=c_2, \la_3(x_3)=x_3$;\qquad\,%
$B_3$: $\la_2(x_2)=c_2, \la_3(x_3)=c_3$.
\item[\bf(C)]
$C_0$: $\la_3(x_3)=x_3$;\qquad\qquad
$C_1$: $\la_3(x_3)=c_3$.
\item[\bf(D)]
$D_0$: $a(x_1,x_2)=x_1, b(x_1,x_2)=x_2, \la_3(x_3)=x_3$;\\
$D_1$: $a(x_1,x_2)=x_1, b(x_1,x_2)=x_2, \la_3(x_3)=c_3$;\\
$D_2$: $a(x_1,x_2)=c_1, b(x_1,x_2)=c_2, \la_3(x_3)=x_3$;\\
$D_3$: $a(x_1,x_2)=c_1, b(x_1,x_2)=c_2, \la_3(x_3)=c_3$.
 \end{enumerate}

We conclude that the equations A, B, C, D are integrable deformations of the equations $A_3$, $B_3$, $C_1$, $D_3$
that will be shown to be most symmetric inside the corresponding family.

\abz\label{defiND}
\begin{defi}\rm
A solution $f$ of any of the equations A, B, C, D on an open set $U\subset M$ with coordinates $(x_1,x_2,x_3)$ is called {\em nondegenerate} if the
corresponding one-form $\al^\la\in T^*U$ defines a Veronese curve at any $x\in U$ (equivalently: the curve
$\la\mapsto\al^\la=\al_0+\la\al_1+\la^2\al_2$ does not lie in any plane, i.e., the 1-forms $\al_0, \al_1, \al_2$ are linearly independent at any point).
\end{defi}

\abz\label{theoP1}
\begin{theo} A generic solution $f$ of any of the equations A, B, C, D is nondegenerate on a small open set $U$. If $f$ is such a solution,
then the corresponding one-form $\al^\la$ defines a Veronese web $\F_\la$ on $U$ by $T\F_\la=\langle\al^\la\rangle^\bot$.
\end{theo}

Here by a generic solution we mean a solution with a generic jet in the Cauchy problem setup.

\medskip

\noindent {\sc Proof.}
The proof of the second statement is the same as that of Theorem \ref{6}(1).
Let us only explain why a generic solution of the equations A-D is nondegenerate.
The condition of degeneracy is vanishing of the determinant of the coefficients of $\al_0,\al_1,\al_2$, so it is a first order PDE that is cubic in the first jets.
A general solution of our second order PDE is a solution to this first order PDE if and only if this first order PDE is an intermediate integral.
But our second order PDE has no intermediate integrals because its symbol is a nondegenerate quadric.
\qed

\abz\label{remnonreg}
\begin{rema}\rm
Let us indicate what happens for degenerate Nijenhuis operators, when two eigenvalues corresponding to different Jordan blocks coincide.
For instance, we can consider real semisimple, but not simple case: $J=\op{diag}(\la(x_1,x_2),\la(x_1,x_2),\nu(x_3))$. The corresponding PDE is
 $$
f_{x_1}f_{x_2x_3}-f_{x_2}f_{x_1x_3}+\frac{\la_{x_1}(x_1,x_2)}{\la(x_1,x_2)-\nu(x_3)}f_{x_2}f_{x_3}-
\frac{\la_{x_2}(x_1,x_2)}{\la(x_1,x_2)-\nu(x_3)}f_{x_1}f_{x_3}=0,
 $$
and its symbol $(f_{x_1}\partial_{x_2}-f_{x_2}\partial_{x_1})\cdot\partial_{x_3}$ is decomposable/degenerate (so it produces neither Veronese curve,
nor conformal Einstein-Weyl structure that we will discuss in Section \ref{X}). Similar situation is with other degenerations, that's why these cases were
rejected from the classification.
\end{rema}

\section{Contact transformations of the equations A, B, C, D}

Let us investigate contact symmetries of the PDEs obtained in the previous section.
The following are the results of straightforward computations in Maple (by the classical method of S.\,Lie).

\abz\label{Prop1}
 \begin{prop}
The contact symmetry algebra of PDE (A) with pairwise different constant functions $\la_i=c_i$,
i.e.\	 that of equation $(A_3)$, is generated by the point symmetries
 \begin{equation}\label{cont1}
h_1(x_1)\partial_{x_1}+h_2(x_2)\partial_{x_2}+h_3(x_3)\partial_{x_3}+h_4(f)\partial_f
 \end{equation}
with arbitrary functions $h_1,h_2,h_3,h_4$ of one argument.
The contact symmetry algebra of PDE (A) with variable $\la_i=\la_i(x_i)$, e.g.\ case $(A_0)$, is generated by the point symmetries
 $$
h(f)\d_f+k_1\cdot\Bigl(\tfrac{1}{\la_1'(x_1)}\d_{x_1}+\tfrac{1}{\la_2'(x_2)}\d_{x_2}+\tfrac{1}{\la_3'(x_3)}\d_{x_3}\Bigr)+
k_2\cdot\Bigl(\tfrac{\la_1(x_1)}{\la_1'(x_1)}\d_{x_1}+\tfrac{\la_2(x_2)}{\la_2'(x_2)}\d_{x_2}+\tfrac{\la_3(x_3)}{\la_3'(x_3)}\d_{x_3}\Bigr)
 $$
with arbitrary two constants $k_1,k_2$ and one function $h$ of one argument.
The pseudogroup $\mathcal{G}_A$ of Lie algebra \eqref{cont1} acts on the class (A) with variable $\la_i$ as a locally transitive transformation pseudogroup.
  \end{prop}

We do not specify contact symmetries of the classes $A_1,A_2$, but the number of arbitrary functions parametrizing them
gradually increases from 1 for $A_0$ to 4 for $A_3$.

\abz\label{Prop2}
 \begin{prop}
The contact symmetry algebra of PDE (B) with pairwise different constant functions $\la_i=c_i$,
i.e.\ that of equation $(B_3)$, is independent of $c_i$ and generated by the point symmetries
 \begin{equation}\label{cont2}
(x_1\cdot h'_1(x_2)+h_2(x_2))\d_{x_1}+h_1(x_2)\d_{x_2}+h_3(x_3)\d_{x_3}+h_4(f)\d_f+k\cdot\d_{x_2}
 \end{equation}
with arbitrary functions $h_1,h_2,h_3,h_4$ of one argument and a constant $k$.
The contact symmetry algebra of PDE (B) with variable $\la_i=\la_i(x_i)$, e.g.\ case $(B_0)$, is generated by the point symmetries
 $$
h_1(x_2)\d_{x_1}+h_2(f)\d_f+k_1\cdot\Bigl(\tfrac{x_1\la_2''(x_2)}{\la_2'(x_2)^2}\d_{x_1}-\tfrac{1}{\la_2'(x_2)}\d_{x_2}-\tfrac{1}{\la_3'(x_3)}\d_{x_3}\Bigr)+
k_2\cdot\Bigl(\tfrac{x_1\la_2(x_2)\la_2''(x_2)}{\la_1'(x_1)^2}\d_{x_1}-\tfrac{\la_2(x_2)}{\la_2'(x_2)}\d_{x_2}-\tfrac{\la_3(x_3)}{\la_3'(x_3)}\d_{x_3}\Bigr)
 $$
with arbitrary two constants $k_1,k_2$ and two functions $h_1,h_2$ of one argument.
The pseudogroup  $\mathcal{G}_B$ of Lie algebra \eqref{cont2} acts on the class (B) with variable $\la_i$ locally transitively.
  \end{prop}

The contact symmetries of the equations $B_1,B_2$ both depend on 3 arbitrary functions.

\abz\label{Prop3}
 \begin{prop}
The contact symmetry algebra of PDE (C) with constant $\la_3=c_3$, i.e.\ that of equation $(C_1)$, is independent of $c_3$
and generated by the point symmetries
 \begin{multline}\label{cont3}
(x_1\cdot h'_1(x_3)+\tfrac12x_2^2h_1''(x_3)+x_2h_2'(x_3)+h_3(x_3))\d_{x_1}+(x_2h'_1(x_3)+h_2(x_3))\d_{x_2}+h_1(x_3)\d_{x_3}\\
+h_4(f)\d_f+k\cdot(x_2\d_{x_2}+2x_3\d_{x_3})
 \end{multline}
with arbitrary functions $h_1,h_2,h_3,h_4$ of one argument and a constant $k$.
The contact symmetry algebra of PDE (C) with variable $\la_3=\la_i(x_3)$, e.g.\ case $(C_0)$, is generated by the point symmetries
 \begin{multline*}
\bigl(x_1h_1(x_3)+h_2(x_3)+\tfrac1{\la'_3(x_3)^2}h'_1(x_3)e^{\la'_3(x_3)x_2}-\tfrac{\la''_3(x_3)}{\la'_3(x_3)^3}h_1(x_3)e^{\la'_3(x_3)x_2}\bigr)\d_{x_1}
+\tfrac1{\la'_3(x_3)}h_1(x_3)\d_{x_2}+h_3(f)\d_f-\\
\Bigl(k_2x_2\tfrac{\la''_3(x_3)}{\la'_3(x_3)^2}+(\sigma x_2\la'''_3(x_3)-k_2\la''_3(x_3))\tfrac1{\la'_3(x_3)^3}
-2\sigma(\la'''_3(x_3)+x_2\la''_3(x_3)^2)\tfrac1{\la'_3(x_3)^4}+5\sigma\tfrac{\la''_3(x_3)^2}{\la'_3(x_3)^5}\Bigr)e^{\la'_3(x_3)x_2}\d_{x_1}\\
+\Bigl(\tfrac{k_2}{\la'_3(x_3)}-\sigma x_2\tfrac{\la''_3(x_3)}{\la'_3(x_3)^2}+\sigma\tfrac{\la''_3(x_3)}{\la'_3(x_3)^3}\Bigr)\d_{x_2}
+\tfrac{\sigma}{\la'_3(x_3)}\d_{x_3};
\qquad \sigma=k_1+k_2\la_3(x_3)
 \end{multline*}
with arbitrary two constants $k_1,k_2$ and three functions $h_1,h_2,h_3$ of one argument.
The pseudogroup  $\mathcal{G}_C$ of Lie algebra \eqref{cont3} acts on the class (C) with variable $\la_i$ locally transitively.
  \end{prop}

\abz\label{Prop4}
 \begin{prop}
The contact symmetry algebra of PDE (D) with pairwise different constant functions $\la_i=c_i$,
i.e.\ that of equation $(D_3)$, is generated by the point symmetries
 \begin{equation}\label{cont4}
h_1(x_1,x_2)\d_{x_1}+h_2(x_1, x_2)\d_{x_2}+h_3(x_3)\d_{x_3}+h_4(f)\d_f
 \end{equation}
with arbitrary smooth functions $h_3,h_4$ of one argument and harmonic duals $h_1,h_2$ (altogether four functions of one argument).
The contact symmetry algebra of PDE (D) with variable $\la_3=\la_3(x_3)$ and complex-analytic $\La(z)=a(x_1,x_2)+ib(x_1,x_2)$,
e.g.\ case $(D_0)$, is generated by the point symmetries
 $$
h(f)\d_f+k_1\cdot\Bigl(\op{Re}\bigl(\tfrac{1}{\La(z)}\d_z\bigr)+\tfrac{1}{\la_3'(x_3)}\d_{x_3}\Bigr)+
k_2\cdot\Bigl(\op{Re}\bigl(\tfrac{\La(z)}{\La'(z)}\d_z\bigr)+\tfrac{\la_3(x_3)}{\la_3'(x_3)}\d_{x_3}\Bigr)
 $$
with arbitrary constants $k_1,k_2$ and function $h$ of one argument.
The pseudogroup $\mathcal{G}_D$ of Lie algebra \eqref{cont4} acts on the class (D) with variable $\La(z),\la_3(x_3)$ locally transitively.
  \end{prop}

The contact symmetries of PDEs $D_1$ and $D_2$ depend on 2 and 3 arbitrary functions respectively.

\abz\label{re-EqGr}
\begin{rema}\rm
It is not difficult to integrate these contact Lie algebras to Lie pseudogroups of local transformations. Then it is apparent
that the equivalences from Remark \ref{remB1} form a subgroup of this pseudogroup. The additional infinite part of the
symmetry pseudogroup comes from the freedom in choice of the first integral of the foliation $\F$: $f\mapsto F(f)$.
\end{rema}

Now we can summarize the computations. The structure equations of the derived symmetry algebras imply the following statement.

\abz\label{theoCont}
\begin{theo} The classes of equations A, B, C and D are pairwise nonequivalent with respect to contact transformations. For any class X
among these its quotient by the corresponding equivalence pseudogroup $\mathcal{G}_X$ has no functional parameters. On generic stratum
the quotient is given by equations $A_0$, $B_0$, $C_0$ and $D_0$ respectively.
\end{theo}

\abz\label{remCont}
\begin{rema}\rm
Note however that there are other equations than $A_i$, $B_j$, $C_k$ and $D_l$ obtained in the quotient. For example, in the class A
we obtain PDE
 $$
(x_2^2-x_3^3)f_{x_1}f_{x_2x_3}+(x_3^3-x_1)f_{x_2}f_{x_3x_1}+(x_1-x_2^2)f_{x_3}f_{x_1x_2}=0
 $$
that is contactly nonequivalent to any of these particular equations in a neighborhood of the origin. The complete list of normal forms is
expressed through the normal forms of functions of one argument.
\end{rema}

\section{B\"{a}cklund  transformations}

By Theorem \ref{theoCont} the PDEs from Section \ref{listPDE} of different type A, B, C or D are not contact equivalent, and some equations
within the same type (for instance, the specifications $A_i, B_j, C_k, D_l$) are also non-equivalent. All these equations are however
equivalent with respect to B\"{a}cklund  transformations, and this also signifies integrability.
This section generalizes the results of I.\ Zakharevich \cite{Z} concerning B\"{a}cklund transformations of the Hirota equation.

\abz\label{12}
\begin{theo}
Let $(\la_1,\la_2,\la_3)$ be a triple of nonzero different functions on $\R^3$ with zero mean, and similarly for $(\La_1,\La_2,\La_3)$.
Assume that $(\la_j/\La_j)_{x_i}=0$ for $i\neq j$. Then the formula
 \begin{equation}\label{e7}
\la_1\La_2f_{x_1}F_{x_2}=\la_2\La_1f_{x_2}F_{x_1},\quad \la_1\La_3f_{x_1}F_{x_3}=\la_3\La_1f_{x_3}F_{x_1}
 \end{equation}
defines the B\"{a}cklund transformation between equation (A) from Section \ref{listPDE} and equation
 \begin{equation}\label{AA}
(\La_2(x_2)-\La_3(x_3))F_{x_1}F_{x_2x_3}+(\La_3(x_3)-\La_1(x_1))F_{x_2}F_{x_3x_1}+(\La_1(x_1)-\La_2(x_2))F_{x_3}F_{x_1x_2}=0
 \end{equation}
(obtained by the substitution $\la\to\La$, $f\to F$ in equation (A)). In other words, for any nondegenerate solution $f$ of PDE (A)
any solution $F$ of system \eqref{e7} gives a nondegenerate solution $F$ of PDE \eqref{AA} and vise versa.
\end{theo}

{\sc Proof.} Put $\phi_i:=\la_i/\La_i$. Then the system (\ref{e7}) is equivalent to the following:
 $$
(F_{x_1},F_{x_2},F_{x_3})=\alpha\cdot(\phi_1f_{x_1},\phi_2f_{x_2},\phi_3f_{x_3}),
 $$
where $\alpha$ is a nonvanishing function. This implies the integrability condition $d\om\wedge\om=0$ for the form
$\om:={\phi_1}f_{x_1}dx_1+{\phi_2}f_{x_2}dx_2+{\phi_3}f_{x_3}dx_3$.
The integrability condition has the form 
$\phi_1(\phi_2-\phi_3)f_{x_1}f_{x_2x_3}+\phi_2(\phi_3-\phi_1)f_{x_2}f_{x_3x_1}+
\phi_3(\phi_1-\phi_2)f_{x_3}f_{x_1x_2}=0$.
On the other hand,
 $$
\phi_1(\phi_2-\phi_3)=\psi\la_1,\ \phi_2(\phi_3-\phi_1)=\psi\la_2,\ \phi_3(\phi_1-\phi_2)=\psi\la_3,
 $$
where we put
 $$
\psi=\frac{\la_2\La_3-\la_3\La_2}{\La_1\La_2\La_3}=\frac{\la_3\La_1-\la_1\La_3}{\La_1\La_2\La_3}=\frac{\la_1\La_2-\la_2\La_1}{\La_1\La_2\La_3},
 $$
and the last equalities are due to $\la_1+\la_2+\la_3=0$, $\La_1+\La_2+\La_3=0$.
This implies that the function $f$ satisfies PDE (A), as claimed. Remark that we did not use the fact that $F$ satisfies (\ref{AA}). \qed

\medskip

The theorem provides B\"{a}cklund transformation between all equations of type (A), in particular between the types $A_0,A_1,A_2,A_3$
(that are contactly non-equaivalent).

Based on the same idea, similar results can be proven also for all other types of equations. The formulation is roughly as follows.
Let $f$ be a nondegenerate solution of a nonlinear second order PDE obtained from a Nijenhuis operator $J$ in $U\subset M^3$
by means of the procedure described in Sections \ref{Hiro}, \ref{listPDE}. Then, if a function $F$ satisfies a certain first order system of linear PDE,
it is a solution of another nonlinear second order PDE that is obtained from the Nijenhuis operator $J^{-1}$ in $U$ by the same procedure.

In other words, the estabilished 1-1 correspondence between Veronese webs and nondegenerate solutions to PDEs of type A, B, C, D
(and even their specifications $A_i,B_j,C_k,D_l$, as shall be proved in Section \ref{Realization}) implies B\"{a}cklund transformations between the
solutions of these PDEs.

\section{Einstein-Weyl structures on the solutions of PDEs}
\label{X}

Now let us construct Einstein-Weyl structure corresponding to solutions of our integrable PDEs.
Recall that this structure consists of a conformal structure $[g]$ and a torsion-free connection
$\nabla$ preserving the conformal class:
 $$
\nabla g=\omega\otimes g.
 $$
The 1-form $\omega$ uniquely encodes the connection $\nabla$ on the 3-dimensional manifold, which is going to
be an arbitrary (graph of) solution of the PDE.

Einstein-Weyl structure is expected to exist due to integrability of this dispersionless PDE. For constant $a_i$, i.e., for the Hirota equation
\eqref{Hirota}, such structure was constructed by Dunajski-Kry\'{n}ński \cite{DK}. We found the corresponding structure for variable $a_i$.

\abz\label{thEW}
 \begin{theo}
The following gives a Weyl structure on a 3D-space with coordinates $(x_1,x_2,x_3)$, parametrized by one function $f=f(x_1,x_2,x_3)$,
where $\la_i=\la_i(x_i)$, $i=1,2,3$.
 \begin{gather*}
g=\frac{(\la_2-\la_3)^2f_{x_1}}{f_{x_2}f_{x_3}}dx_1^2 +\frac{(\la_1-\la_3)^2f_{x_2}}{f_{x_1}f_{x_3}}dx_2^2 +
\frac{(\la_1-\la_2)^2f_{x_3}}{f_{x_1}f_{x_2}}dx_3^2+\\
\frac{2(\la_1-\la_3)(\la_2-\la_3)}{f_{x_3}}dx_1dx_2 +\frac{2(\la_1-\la_2)(\la_3-\la_2)}{f_{x_2}}dx_1dx_3
+\frac{2(\la_2-\la_1)(\la_3-\la_1)}{f_{x_1}}dx_2dx_3,\\
\omega=\left(\Bigl(\frac1{\la_1-\la_2}+\frac1{\la_1-\la_3}\Bigr)\la_1'-\Bigl(\frac1{\la_1-\la_2}\frac{f_{x_1}}{f_{x_2}}\Bigr)\la_2'
-\Bigl(\frac1{\la_1-\la_3}\frac{f_{x_1}}{f_{x_3}}\Bigr)\la_3'-\frac{f_{x_1x_1}}{f_{x_1}}\right)dx_1 \\
+\left(\Bigl(\frac1{\la_2-\la_1}+\frac1{\la_2-\la_3}\Bigr)\la_2'-\Bigl(\frac1{\la_2-\la_1}\frac{f_{x_2}}{f_{x_1}}\Bigr)\la_1'
-\Bigl(\frac1{\la_2-\la_3}\frac{f_{x_2}}{f_{x_3}}\Bigr)\la_3'-\frac{f_{x_2x_2}}{f_{x_2}}\right)dx_2\\
+\left(\Bigl(\frac1{\la_3-\la_1}+\frac1{\la_3-\la_2}\Bigr)\la_3'-\Bigl(\frac1{\la_3-\la_1}\frac{f_{x_3}}{f_{x_1}}\Bigr)\la_1'
-\Bigl(\frac1{\la_3-\la_2}\frac{f_{x_3}}{f_{x_2}}\Bigr)\la_2'-\frac{f_{x_3x_3}}{f_{x_3}}\right)dx_3.
 \end{gather*}
It is Einstein-Weyl iff the function $f$ satisfies PDE (A).
 \end{theo}

\abz\label{remEW}
 \begin{rema}\rm
The Lax pair does not contain the derivative by the spectral parameter $\lambda$, so it is an integrable hyper-CR
equation. Thus the procedure of Jones-Tod \cite{JT} gives a way to construct Einstein-Weyl structure with the trick described
in \cite{DK}. We omit detils of this computation.

It is interesting to note that the same Einstein-Weyl structure is obtained via the "universal" formula of \cite{FK}:
 $$
\omega_k=2g_{kj}D_{x_l}(g^{lj})+D_{x_k}\log\det(g_{ij}).
 $$
Here $D_{x_l}$ are the total derivative operators \cite{KLV,KL}.

Notice that equation (A) is not translationary invariant, which is the standard setup for application of the method of hydrodynamic
integrability. In particular, the results of \cite{FK} formally do not apply to this equation. Yet the formula (miraculously) works here as well.
In a similar way we derive Einstein-Weyl structures for integable PDE B, C and D (the formula for $\omega_k$ works
for case D).
 \end{rema}

Existence of Einstein-Weyl structures on the solutions of equations A, B, C and D exhibits them as reductions of the universal Einstein-Weyl
equation, which is integrable by the twistor methods \cite{H,DFK}. This again confirms integrability of our deformations. The Einstein-Weyl structures
could be computed explicitely for all types, similarly as it is done in case (A). Yet the existence of this structure on the solutions
of these equations follow from the general result of \cite{CK}, since our Lax pairs are easily checked to be characteristic, i.e., null for the
canonical conformal structure.

\section{ Realization theorem}\label{Realization}

The aim of this section is to prove analogs of Theorem \ref{6}(2) for A, B, C, D equations.
We begin with the general situation, and then specify to the 3-dimensional case.

 \abz\label{defSP}
\begin{defi}\rm
Consider a Veronese web $\F_\la$ on a manifold $M^{n+1}$, given by $T\F_\la=\langle \al^\la\rangle^\bot$,
where $\al^\la=\al_0+\la\al_1+\cdots+\la^n\al_n$ and $\al_0,\al_1,\ldots,\al_n$ is a local coframe on an open set $U \subset M$.
A smooth function $\phi:U\to \R$ is called {\em self-propelled}
if $d\phi$ is proportional to $\al^\phi$. If the coefficient of proportionality is nonzero, we denote this by $d\phi\sim\al^\phi$.
But the coefficient is allowed to be zero, so a constant function is also considered self-propelled.
 \end{defi}

\abz\label{lemmSP}
\begin{lemm}
Let $\F_\la$ be a Veronese web on $M^{n+1}$. Then in a vicinity of any point $x\in M$ there exist $n+1$ functionally independent self-propelled
functions $\phi_0(x),\phi_1(x),\dots,\phi_n(x)$. If $X_0,\dots,X_n$ is the frame dual to the coframe $\al_0,\dots,\al_n$ defining the Veronese web,
the condition on the function $\phi$ to be self-propelled is the following system PDEs:
 \begin{equation}\label{sys}
\phi X_0\phi=X_1\phi, \ \dots,\ \phi X_{n-1}\phi=X_n\phi.
 \end{equation}
\end{lemm}

\noindent {\sc Proof.} The required relation $\al_0+\dots+\phi^n\al_n\sim (X_0\phi) \al_0+\dots+(X_n\phi) \al_n$ is equivalent to
vanishing of the determinants
$$
\left|
  \begin{array}{cc}
  1 & \phi \\ X_0\phi  & X_1\phi  \\
  \end{array}
\right|,\ \dots,\
\left|
  \begin{array}{cc}
  \phi^{n-1} & \phi^n \\ X_{n-1}\phi  & X_n\phi  \\
  \end{array}
\right|,
$$
which is equivalent to system (\ref{sys}). Let $F(x_1,\dots,\la)$ be a $\la$-parametric first integral of the folitation $\F_\la$.
The following formula gives a family of implicit solutions $\phi(x)$ of system \eqref{sys} depending on an arbitrary smooth function of one
variable $f=f(\la)$ that locally satisfies $f'(\la)\not=F_\la$.
 \begin{equation}\label{eqt}
F(x,\phi(x))=f(\phi(x)).
 \end{equation}
Indeed, differentiating this equality along $X_k-\phi(x) X_{k-1}$ and writing $x=(x_1,\dots,x_n)$ we get
 \begin{equation}\label{eqt1}
d_xF(x,\lambda)(X_k-\la X_{k-1})|_{\lambda=\phi(x)}+(F_\la(x,\phi(x))-f'(\phi(x)))\cdot(X_k\phi(x)-\phi(x) X_{k-1}\phi(x))=0.
 \end{equation}
The first term vanishes since $X_k-\la X_{k-1}\in\langle\al^\la\rangle^\perp$, and the claim follows.

Choosing $n$ solutions $\phi_0,\ \dots,\ \phi_n$ with initial values $c_0,\dots,c_n$ at $x\in M$ being
pairwise different and with nonzero $\psi_i:=X_0\phi_i|_x$, we compute from (\ref{sys}) the Jacobian at $x$:
 $$
\mathrm{Jac}_{x}(\phi_0,\phi_1,\dots,\phi_n)\sim
\left|\begin{array}{cccc}\psi_0 & c_0\psi_0 & \dots & c_0^n\psi_0\\ \psi_1 & c_1\psi_1 & \dots & c_1^n\psi_1\\
\vdots & \vdots & \ddots & \vdots \\
\psi_n & c_n\psi_n & \dots & c_n^n\psi_n\end{array}\right|
=\psi_0\psi_1\cdots\psi_n\left|
\begin{array}{cccc}1 & c_0 & \dots & c_0^n\\ 1 & c_1 & \dots & c_1^n\\
\vdots & \vdots & \ddots & \vdots\\ 1 & c_n & \dots & c_n^n\end{array}\right|.
 $$
Since the Vandermonde determinant with the second column consisting of pairwise different entries is nonzero,
we obtain $n$ functionally independent solutions of (\ref{sys}). \qed

\abz\label{remmk}
\begin{rema}\rm
Below we will also need a generalization of the notion of a self-propelled function (we specify to the case $n=2$). A smooth complex valued function $\phi(x)=\eta(x)+i\zeta(x)$  on $M$ is {\em self-propelled} if it satisfies system \ref{sys}, or, equivalently if the real-valued functions $\eta(x)$ and $\zeta(x)$ satisfy the following system of equations
\begin{eqnarray}\label{sys1}
\eta X_0 \eta-\zeta X_0\zeta=X_1\eta, \eta X_0 \zeta+\zeta X_0\eta=X_1\zeta,\\
\eta X_1 \eta-\zeta X_1\zeta=X_2\eta, \eta X_1 \zeta+\zeta X_1\eta=X_2\zeta.\nonumber
\end{eqnarray}
We will also say that a {\em pair} of functions $(\eta,\zeta)$ is {\em self-propelled} if it satisfies the system above.

Motivated by the case D of 3D Nijenhuis operators, which is necessary analytic in $z=x_1+ix_2$,
we will assume real-analyticity for the PNO data of the corresponding case.
Working in real-analytic category one can prove the local existence of complex self-propelled functions using the same arguments as in the real case
for the complexification of equation (\ref{eqt}). Similarly, we can show that there exist local functionally independent real(-analytic) functions
$\eta,\zeta,\psi$ such that both the pair functions $(\eta,\zeta)$ and the function $\psi$ are self-propelled.
\end{rema}

\abz\label{lemBKlie}
\begin{lemm}
Let $c_{ij}^k$ be the structure functions of the frame $X_0,\dots,X_n$: $[X_i,X_j]=c_{ij}^kX_k$.
Then the compatibility conditions (= necesary and sufficient conditions of solvability with any admissible Cauchy data) of
system \eqref{sys} is the vanishing of coefficients of the polynomials
 \begin{equation}\label{sysCijk}
\Psi_{ijk}(\phi)=\sum_{m=0}^n\bigl(c_{ij}^m\phi^{k+m}+c_{jk}^m\phi^{i+m}+c_{ki}^m\phi^{j+m}\bigr)
 \end{equation}
for all triple of numbers $(i,j,k)$ from $[0..n]$.
\end{lemm}

\noindent {\sc Proof.} We will use the language of the geometry of PDE \cite{S,KLV}.
The symbol of system \eqref{sys} is (in this proof $T=T_xM$ is the tangent space pulled-back to the point of the equation)
 $$
g_1(\phi)=\langle X_1-\phi X_0,\dots,X_n-\phi X_{n-1}\rangle^\perp =\langle d\phi:d\phi([X_0,\dots,X_n])=t\cdot[1,\phi,\dots,\phi^n]\rangle\subset T^*.
 $$
This space is one-dimensional and its prolongations $g_k=g_1\otimes S^{k-1}T^*\cap S^kT^*$
are one-dimensional as well. In fact, $g_k=\langle (\al^\phi)^k\rangle$, $\al^\phi\in T^*$. The Spencer $\delta$-sequence
 $$
\dots\stackrel{\delta}\to g_{i+1}\otimes \Lambda^{j-1}T^*\stackrel{\delta}\to g_{i}\otimes \Lambda^{j}T^*
\stackrel{\delta}\to g_{i-1}\otimes \Lambda^{j+1}T^* \stackrel{\delta}\to \dots
 $$
has cohomology $H^{i,j}$ at the term $g_{i}\otimes \Lambda^{j}T^*$ \cite{S,KL}.
The only non-trivial second cohomology groups (encodes the compatibility) is
$H^{0,2}=\op{Ker}(\d_\phi:\Lambda^2T^*\to\La^3T^*,\omega\mapsto\omega\wedge\al^\phi)$,
they live on 2-jets.

These compatibility conditions visualise as follows.
The first prolongation of system \eqref{sys} written as $\{X_k\phi=\phi^k\sigma\}$ implies
$\kappa_{ij}:=[X_i,X_j]\phi+(i-j)\phi^{i+j-1}\sigma^2=(\phi^jX_i-\phi^iX_j)\cdot\sigma$. The cocycle condition
$\phi^l\kappa_{ij}+\phi^i\kappa_{jl}+\phi^j\kappa_{li}=0$ implies compatibility conditions \eqref{sysCijk}.

This implies in general, by the Cartan-K\"ahler theory \cite{KLV}, local integrability only provided the data are analytic.
However our system has formal solution depending on 1 function of 1 variable, and hence
here we can exploit a Sophus Lie theorem \cite{L,De}, which implies that in this case there exists a local solution
in the smooth category, see \cite{K}. \qed

\bigskip

In the case of our interest $\dim M=3$ ($n=2$) there is only one polynomial $\Psi_{012}$ and we conclude with respect to the frame $X_0,X_1,X_2$:
 $$
c_{12}^0+(c_{12}^1-c_{02}^0)\phi+(c_{01}^0-c_{02}^1+c_{12}^2)\phi^2+(c_{01}^1-c_{02}^2)\phi^3+c_{01}^2\phi^4=0.
 $$
Vanishing of this polynomial in $\phi$ is equivalent to such structure relations:
 \begin{gather}
[X_0,X_1]=b_0X_0+b_1X_1,\ [X_1,X_2]=c_1X_1+c_2X_2, \label{str1}\\
 [X_0,X_2]=c_1X_0+(c_2+b_0)X_1+b_1X_2. \label{str2}
 \end{gather}

\abz\label{mainTh}
\begin{theo}
Let $(\F,\bar{J})$ be a partial Nijenhuis operator of generic type (see Theorem \ref{theo1}) on a 3-dimensional manifold $M$. Then in a
neighborhood $U$ of every point $x\in M$ there exists a Nijenhuis operator $J:TM\to TM$ of any type A, B, C or D (in the last case the PNO
is assumed real-analytic) such that $J|_{T\F}=\bar{J}$.
\end{theo}

\noindent {\sc Proof.} Consider $(\F,\bar{J})$ in $U$. The intersection $D_1:=T\F\cap \bar{J}T\F$ is a one-dimensional distribution.
Choose a nonvanishing vector field $X_1\in\Gamma(D_1)$ and put $X_0:=\bar{J}^{-1}X_1,X_2:=\bar{J}X_1$.
Then $X_0,X_1,X_2$ is a frame satisfying the structure equations \eqref{str1}-\eqref{str2} for some functions $b_0,b_1,c_1,c_2$.

The first line \eqref{str1} is due to the integrability of the distributions $T\F$ and $NT\F$.
To prove \eqref{str2} decompose $[X_0,X_2]=d_0X_0+d_1X_1+d_2X_2$ and use the definition of a PNO:
by condition 1 of this definition we have $[X_0,X_1]_{\bar{J}}=d_0X_0+(d_1-b_0)X_1+(d_2-b_1)X_2\in T\F$, which implies $d_2=b_1$;
by condition 2 of this definition we have $c_1X_1+c_2X_2=[\bar{J}X_0,\bar{J}X_1]=\bar{J}([X_0,X_1]_{\bar{J}})=d_0X_1+(d_1-b_0)X_2$,
which implies $d_0=c_1, d_1=c_2+b_0$.

The matrix of the operator $\bar{J}:T\F\to TM$ with respect to the bases $X_0,X_1$ in $T\F$ and $X_0,X_1,X_2$ in $TM$ is equal to
 $$
\left[
  \begin{array}{cc}
    0 &   0  \\     1 & 0  \\     0 & 1  \\
  \end{array}
\right].
 $$
Define $J$ by $J|_{T\F}=\bar{J}$ and $JX_2=f_0 X_0+f_1X_1+f_2X_2$, where $f_i$ are smooth functions on $U$.
Thus the matrix of $J$ in the frame $X_0,X_1,X_2$ is
 $$
\left[
  \begin{array}{ccc}
    0 & 0 & f_0 \\     1 & 0 & f_1 \\     0 & 1 & f_2 \\
  \end{array}
\right].
 $$
Direct calculations show that $N_{J}(X_1,X_2)=0$ is equivalent to
 \begin{equation}\label{er1}
X_2f_0=f_0X_1f_2, \quad X_2f_1=X_1f_0+f_1X_1f_2, \quad X_2f_2=X_1f_1+f_2X_1f_2,
 \end{equation}
and, analogously, $N_{J}(X_0,X_2)=0$ is equivalent to
 \begin{equation}\label{er2}
X_1f_0=f_0X_0f_2, \quad X_1f_1=X_0f_0+f_1X_0f_2, \quad X_1f_2=X_0f_1+f_2X_0f_2.
 \end{equation}
Now let $f_0=\phi_1\phi_2\phi_3$, $f_1=-(\phi_1\phi_2+\phi_1\phi_3+\phi_2\phi_3)$, $f_3=\phi_1+\phi_2+\phi_3$ for some
functions $\phi_1,\phi_2,\phi_3$. Then it is easy to see that once the functions $\phi_i$ satisfy the system of equations (\ref{sys}),
the functions $f_i$ satisfy the systems of equations (\ref{er1}), \ref{er2}). In other words, if the functions $\phi_1,\phi_2,\phi_3$ are self-propelled, the Nijenhuis tensor $N_{J}$ of the (1,1)-tensor $J$ given in the frame $X_0,X_1,X_2$ by the matrix
\begin{equation}\label{matr}
F(\phi_1,\phi_2,\phi_3):=\left[
  \begin{array}{ccc}
  0 & 0 &  \phi_1\phi_2\phi_3\\
    1 & 0 & -\phi_1\phi_2-\phi_1\phi_3-\phi_2\phi_3 \\
    0 & 1 & \phi_1+\phi_2+\phi_3 \\
  \end{array}
\right]
\end{equation}
vanishes (recall that $N_{J}(X_0,X_1)=N_{\bar{J}}(X_0,X_1)=0$ by the assumptions of the theorem).

Now if we take $(\phi_1,\phi_2,\phi_3)$ to be $(\la_1,\la_2,\la_3)$, $(\la_2,\la_2,\la_3)$ or
$(\la_3,\la_3,\la_3)$, where $\la_i$ are functionally independent and self-propelled (and we can write $\la_i=\la_i(x_i)$),
we obtain $J$ of type A, B, or C respectively. To get type D use Remark \ref{remmk}. \qed

\abz\label{A_iB_jC_kD_l}
\begin{rema}\rm
We can even get more specified forms $A_i$, $B_j$, $C_k$ for the first three types (in notations of Section \ref{listPDE}). For this
let $\la_i=\la_i(x_i)$ be functionally independent self-propelled functions, $c_i$ pairwise different constants and define $J$ as $F_X$,
where $X$ denotes the type:
 \begin{itemize}
\item[(A)] $F_{A0}:=F(\la_1,\la_2,\la_3)$; $F_{A1}:=F(\la_1,\la_2,c_3)$; $F_{A2}:=F(\la_1,c_2,c_3)$; $F_{A3}:=F(c_1,c_2,c_3)$;
\item[(B)] $F_{B0}:=F(\la_2,\la_2,\la_3)$; $F_{B1}:=F(\la_2,\la_2,c_3)$; $F_{B2}:=F(c_2,c_2,\la_3)$; $F_{B3}:=F(c_2,c_2,c_3)$;
\item[(C)] $F_{C0}:=F(\la_3,\la_3,\la_3)$; $F_{C1}:=F(c_3,c_3,c_3)$.
 \end{itemize}
In the case of real-analytic PNO, we can also realize it by the type $D_l$.

Taking $a=a(x_1,x_2)$, $b=b(x_,x_2)$, $\la_3=\la_3(x_3)$ functionally independent functions such that $(a,b)$ is a self-propelled pair
and $\la_3$ is self-propelled itself (see Remark \ref{remmk}) and constants $c_1,c_2\not=0,c_3$ one can put also
 \begin{itemize}
\item[(D)]  $F_{D0}:=F(a+ib,a-ib,\la_3)$; $F_{D1}:=F(a+ib,a-ib,c_3)$;
\item[]  $F_{D2}:=F(c_1+ic_2,c_1-ic_2,\la_3)$; $F_{D3}:=F(c_1+ic_2,c_1-ic_2,c_3)$.
 \end{itemize}
The matrices $F_X$ are the Frobenius forms of all the Nijenhuis operators listed in Appendix.
\end{rema}

\abz\label{corRR}
\begin{coro}
Let $\F_\la$ be a Veronese web on a 3-dimensional smooth manifold $M$. Then for any type $(X)$ of the equations listed in Section \ref{listPDE} in a neighbourhood of any point on $M$ there exist local coordinates $(x_1,x_2,x_3)$ such that any first integral $f$ of the foliation $\F_\infty$
expressed in these coordinates is a nondegenerate solution of equation $(X)$ (see Definition \ref{defiND}) for $X=A_i,B_j,C_k\text{ or }D_l$.
\end{coro}

\noindent {\sc Proof.} Let $\F_\la$ be a Veronese web and let $f$ be a function such that $\F_\infty=\{f=const\}$. Consider a PNO
$\bar{J}:T\F_\infty\to TM$ with $\op{Im}\bar{J}= T\F_0$ which corresponds to $\F_\la$ by Theorem \ref{theo1}.
Repeat the construction from the proof of Theorem \ref{mainTh} to get a Nijenhuis operator $J$, $J|_{T\F_\infty}=\bar{J}$, and the coordinates $(x_1,x_2,x_3)$ such that the matrix of $J$ in the basis $\{\partial_{x_i}\}$ has the form $J_X$ from the list of Appendix. The distribution $J(T\F_\infty)=T\F_0$ is integrable, hence $d\om\wedge\om=0$ and $f$ expressed in coordinates $(x_1,x_2,x_3)$ is a solution of the corresponding equation $(X)$. Nondegeneracy of $f$ follows from Theorem \ref{theo1}, since the degeneracy would imply that $J|_{T\F_\infty}$ is of nongeneric type. \qed

\section{Examples}

We want to illustrate the relation between Veronese webs and PDE, and show how this can be used to construct exact solutions.
In this we will be following the proofs of Lemma \ref{lemmSP} and  Theorem \ref{mainTh}.

As an example we consider $\mathfrak{sl}_2(\R)$ represented by a vector frame $X_0,X_1,X_2$ with commutation relations
$[X_0,X_1]=X_0, [X_1,X_2]=X_2, [X_0,X_2]=2X_1$. Note that these vector fields satisfy equations \eqref{str1}-\eqref{str2}
crusial for Theorem (\ref{mainTh}). Denoting by $\al_0,\al_1,\al_2$ is the dual coframe to $X_0,X_1,X_2$, we conclude that
the 1-form $\al^\la:=\al_0+\la\al_1+\la^2\al_2$ is integrable and so defines a Veronese web.

This Veronese web is nonflat (i.e., in no coordinate system the leaves of the foliations $\F_\la$ are parallel planes) due to the nonintegrability
of the distribution $\langle X_0,X_2\rangle$. Let us choose the following realization of the frame in $\R^3$ (away from $0$):
$$
X_0=\partial_{x_1}+\partial_{x_2}+\partial_{x_3},\
X_1=x_1\partial_{x_1}+x_2\partial_{x_2}+x_3\partial_{x_3},\
X_2=x_1^2\partial_{x_1}+x_2^2\partial_{x_2}+x_3^2\partial_{x_3}.
$$

\abz\label{exaa1}
\begin{exa}\rm
The function $F(x_1,x_2,x_3,\la)$, mentioned in the proof of Lemma \ref{lemmSP}, whose level sets coincide with the leaves of $\F_\la$,
i.e., the solution of the system of equations
 $$
(X_1-\la X_0)F(x_1,x_2,x_3,\la)=0,\ (X_2-\la X_1)F(x_1,x_2,x_3,\la)=0
 $$
is given by the "cross-ratio" formula
 $$
F(x_1,x_2,x_3,\la)=\frac{(x_3-x_2)(x_1-\la)}{(x_1-x_2)(x_3-\la)}.
 $$
From this we get the implicit solution of system \eqref{sys}: $F(x_1,x_2,x_3,\phi(x_1,x_2,x_3))=f(\phi(x_1,x_2,x_3))$.
Taking $f=c=\op{const}$ we obtain the following 1-parametric family of the explicit solutions of this system, i.e.,
of self-propelled functions for the initial Veronese web:
 $$
\phi(x_1,x_2,x_3)=\frac{cx_3(x_1-x_2)+x_1(x_3-x_2)}{c(x_1-x_2)+(x_3-x_2)}.
 $$
In particular, the solutions corresponding to $c=0,-1,\infty$,
 $$
\psi_1=x_1,\ \psi_2=x_2,\ \psi_3=x_3,
 $$
are the original coordinates. The function $F(x_1,x_2,x_3,\infty)=\frac{x_3-x_2}{x_1-x_2}$ ``cutting'' the foliation $\F_\infty$
is a particular solution of the equation of type $(A_0)$
 $$
(x_2-x_3)f_{x_1}f_{x_2x_3}+(x_3-x_1)f_{x_2}f_{x_3x_1}+(x_1-x_2)f_{x_3}f_{x_1x_2}=0.
 $$
Now fix $\la=c$ and take $x_1,x_2,\bar{x}_3=F(x_1,x_2,x_3,c)$ as coordinates. Then
$F(x_1,x_2,x_3,\infty)=\frac{(c-x_2)\bar{x}_3}{(x_1-x_2)\bar{x}_3+(x_1-c)}$ is a solution of the equation of type $(A_1)$
 $$
(x_2-c)f_{x_1}f_{x_2\bar{x}_3}+(c-x_1)f_{x_2}f_{\bar{x}_3x_1}+(x_1-x_2)f_{\bar{x}_3}f_{x_1x_2}=0.
 $$
Analogously one can build solutions of equations of types $(A_2)$, $(A_3)$.

\end{exa}

\abz\label{exaa2}
\begin{exa}\rm

To construct a solution of equation $(D_0)$ we will first build a self-propelled pair of functions $(\psi_1,\psi_2)$ (see Remark \ref{remmk}).
To this end solve the equation $F(x_1,x_2,x_3,\psi_1(x_1,x_2,x_3)+i\psi_2(x_1,x_2,x_3))=c_1+ic_2$ with respect to $\psi_1,\psi_2$
for some real constants $c_1,c_2$, say $c_1=c_2=1$:
 $$
\psi_1=\frac{(x_1^2+x_2x_3)(x_2+x_3)-4x_1x_2x_3}{(x_1-x_2)^2+(x_1-x_3)^2},\
\psi_2=\frac{(x_2-x_3)(x_1-x_3)(x_1-x_2)}{(x_1-x_2)^2+(x_1-x_3)^2}.
 $$
Now take the coordinates $\bar{x}_1=\psi_1$, $\bar{x}_2=\psi_2$, $x_3$ and express $F(x_1,x_2,x_3,\infty)$ in them:
$\frac{x_3-x_2}{x_1-x_2}=f(\bar{x}_1,\bar{x}_2,x_3)=\frac{\bar{x}_1-\bar{x}_2-x_3}{\bar{x}_2}$. This is a solution of equation $(D_0)$ with
substituted $(x_1,x_2,x_3)=(\bar{x}_1,\bar{x}_2,x_3)$.

\smallskip

In general, in order to find a solution of equation $(X)$, with $X=A_i, B_j, C_k\text{ or }D_l$ one should first of all find coordinates
$(\bar{x}_1,\bar{x}_2,\bar{x}_3)$ in which the Nijenhuis operator $J$ built in Theorem \ref{mainTh} takes the canonical form.
Then the function $F(x_1,x_2,x_3,\infty)=f(\bar{x}_1,\bar{x}_2,\bar{x}_3)$ expressed in these coordinates will be a solution of this equation.
The coordinates can be defined intrinsically and uniquely with respect to the operator $J$,
for instance, being the eigenfunctions of $J$. In such a case it is easy to find them. However, in some of the cases, especially that containing
higher-dimensional Jordan blocks, they are defined nonuniquely and one needs more efforts to find them.

\end{exa}

\abz\label{exaa3}
\begin{exa}\rm

Consider case $(B_0)$. Here the coordinates $x_2,x_3$ are the eigenfunctions of $J$, which can be any functionally independent self-propelled functions.
The corresponding eigenvectors can be calculated by the formula $Z_i=P_i^0X_0+P_i^1X_1+P_i^2X_2$, $i=1,3$, where $P_i^j$
are the entries of the matrix $P_{B0}=P_{B0}(x_2,x_3)$ (see the Appendix), which gives
 \begin{eqnarray*}
Z_1=\frac{(x_1-x_2)(x_1-x_3)}{x_2-x_3}\cdot\partial_{x_1},\
Z_3=\frac{(x_1-x_2)^2}{(x_2-x_3)^2}\cdot\partial_{x_1}+\partial_{x_3}.
 \end{eqnarray*}
These vector fields should be proportional to $\partial_{\bar{x_1}},\partial_{\bar{x_3}}$ respectively and $\partial_{\bar{x_2}}$ should be
adjoint to $\partial_{\bar{x_1}}$, i.e., $(J-x_2\Id)\partial_{\bar{x_2}}=\partial_{\bar{x_1}}$. On the other hand, the vector field
 $$
Z_2=P_2^0X_0+P_2^1X_1+P_2^2X_2=\frac{-2x_2x_3+x_3^2+2x_1x_2-x_1^2}{(x_2-x_3)^2}\cdot\partial_{x_1}+\partial_{x_2}
 $$
is adjoint to $Z_1$ and the system of equations
$$
Z_i\bar{x}_j=\delta_{ij}
$$
has a unique (up to additive constants) solution $\bar{x}_1=\ln\frac{(x_1-x_2)(x_2-x_3)}{x_1-x_3}$, $\bar{x}_2=x_2$, $\bar{x}_3=x_3$,
hence we have found  needed coordinates and $Z_i=\partial_{\bar{x_i}}$. The function
$F(x_1,x_2,x_3,\infty)=f(\bar{x}_1,\bar{x}_2,\bar{x}_3)=1-\bar{x}_2e^{-\bar{x}_1}+\bar{x}_3e^{-\bar{x}_1}$ is a solution of equation $(B_0)$.

\end{exa}

\abz\label{exaa4}
\begin{exa}\rm

Consider case $(C_0)$. Here the coordinate $x_3$ is an eigenfunction, which can be any self-propelled function. Set $\bar{x}_3=x_3$.
The basis fields $\partial_{\bar{x}_i}$ in the coordinate system we are looking for satisfy the following relations
$(J-\bar{x}_3\Id)\partial_{\bar{x}_2}=0$, $(J-\bar{x}_3\Id)\partial_{\bar{x}_1}=\partial_{\bar{x}_2}$,
$(J-\bar{x}_3\Id)\partial_{\bar{x}_3}=\partial_{\bar{x}_1}-\bar{x}_2\partial_{\bar{x}_2}$.
The vector fields $Z_i=P_i^0X_0+P_i^1X_1+P_i^2X_2$, $i=1,2,3$, where $P_i^j$ are the entries of the matrix $P_{C0}=P_{C0}(\bar{x}_2,x_3)$
(see the Appendix) obey the same relations: $(J-\bar{x}_3\Id)Z_2=0$, $(J-\bar{x}_3\Id)Z_1=Z_2$, $(J-\bar{x}_3\Id)Z_3=Z_1-\bar{x}_2Z_2$;
here $\bar{x}_2=g(x_1,x_2,x_3)$ is an unknown function. Explicitly,
 \begin{gather*}
Z_1=(x_1-x_3)\partial_{x_1}+(x_2-x_3)\partial_{x_2},\
Z_2=(x_1-x_3)^2\partial_{x_1}+(x_2-x_3)^2\partial_{x_2},\\
Z_3=-(g(x_1-x_3)-1)\partial_{x_1}-(g(x_2-x_3)-1)\partial_{x_2}+\partial_{x_3}.
 \end{gather*}
We are not so lucky as in the previous case, since these vector fields do not pairwise commute and depend on an unknown function.
The following vector fields $Z_1'=Z_1-gZ_2$, $Z_2'=Z_2$, $Z_3'=Z_3+gZ_1$, where $g=1/(x_1-x_3)$ or, explicitly,
 $$
Z_1'=\frac{(x_1-x_2)(x_2-x_3)}{x_1-x_3}\cdot\partial_{x_2},\
Z_2'=(x_1-x_3)^2\partial_{x_1}+(x_2-x_3)^2\partial_{x_2},\
Z_3'=\partial_{x_1}+ \partial_{x_2} + \partial_{x_3},
 $$
obey the same relations and pairwise commute. The solution of the system $Z_i'\bar{x}_j=\delta_{ij}$ is
 $$
\bar{x}_1=\ln\frac{(x_1-x_3)(x_2-x_3)}{x_1-x_2},\ \bar{x}_2=g=\frac{1}{x_1-x_3},\ \bar{x}_3=x_3.
 $$
The function $F(x_1,x_2,x_3,\infty)$ is expressed in these new coordinates is equal to $\bar{x}_2e^{\bar{x}_1}$,
which gives a nondegenerate solution of equation $(C_0)$.

\end{exa}

\section{Appendix: Classification of cyclic Nijenhuis operators in 3D (after F.\,J.\ Turiel)}\label{app}

In papers \cite{T,GM} the authors obtained a local classification of Nijenhuis operators $J:TM\to TM$ (in a vicinity of a regular point \cite[p. 451]{T})
under the additional assumption of existence of a complete family of the conservation laws. This assumption is equivalent to vanishing of
the invariant $P_J$, which is automatically trivial in the case of cyclic $J$ \cite[p. 450]{T}, i.e., when the space $T_xM$ is cyclic for $J_x$, for any $x\in M$.
Here we recall the normal forms obtained in this case for 3-dimensional $M$.
As we stated in Section \ref{3DNO} no additional assumption (like "cyclic") is needed in 3D to conclude these forms.

The results of \cite{T} imply that for any (cyclic) Nijenhuis operator ($(1,1)$-tensor) in a vicinity of a regular/generic point $x^0$ there exists a
local system of coordinates $(x_1,x_2,x_3)$ centered around $(x^0_1,x^0_2,x^0_3)$ and pairwise distinct constants $c_1,c_2,c_3$
($c_i\neq x_i^0$), such that the matrix of $J$ in the basis $\{\partial_{x_i}\}$ is one from the following list, where we also indicate
the cyclic Frobenius forms $F$ as well as the operators $P$ for which $PJP^{-1}=F$.
\begin{itemize}
\item[\bf ${\bold A_0}$.] $J_{A0}=J_{A0}(x_1,x_2,x_3):=\left[
           \begin{array}{ccc}
             x_1  & 0 & 0 \\
             0 & x_2 & 0 \\
             0 & 0 & x_3 \\
           \end{array}
         \right]$,\\
$F_{A0}=F_{A0}(x_1,x_2,x_3)=\left[
           \begin{array}{ccc}
             0 & 0 & x_1x_2x_3 \\
             1 & 0 & -x_1x_2-x_1x_3-x_2x_3 \\
             0 & 1 & x_1+x_2+x_3 \\
           \end{array}
         \right]$,\\
$P_{A0}=P_{A0}(x_1,x_2,x_3)=\left[\begin{array}{ccc}
{\frac{x_2x_3}{(x_1-x_2)(x_1-x_3)}} & {\frac{x_1x_3}{(x_2-x_1)(x_2-x_3)}} & {\frac{x_1x_2}{(x_3-x_1)(x_3-x_2)}}\\
\noalign{\medskip}{\frac{-(x_2+x_3)}{(x_1-x_2)(x_1-x_3)}} & {\frac{-(x_1+x_3)}{(x_2-x_1)(x_2-x_3)}} & {\frac{-(x_1+x_2)}{(x_3-x_1)(x_3-x_2)}}\\
\noalign{\medskip}{\frac{1}{(x_1-x_2)(x_1-x_3)}} & {\frac{1}{(x_2-x_1)(x_2-x_3)}} & {\frac{1}{(x_3-x_1)(x_3-x_2)}}
 \end{array}\right]$.
\item[\bf ${\bold A_1}$.] $J_{A1}:=J_{A0}(x_1,x_2,c_3)$, $F_{A1}:=F_{A0}(x_1,x_2,c_3)$, $P_{A1}=P_{A0}(x_1,x_2,c_3)$.
\item[\bf ${\bold A_2}$.] $J_{A2}:=J_{A0}(x_1,c_2,c_3)$, $F_{A2}:=F_{A0}(x_1,c_2,c_3)$, $P_{A2}=P_{A0}(x_1,c_2,c_3)$.
\item[\bf ${\bold A_3}$.] $J_{A3}:=J_{A0}(c_1,c_2,c_3)$, $F_{A3}:=F_{A0}(c_1,c_2,c_3)$, $P_{A3}=P_{A0}(c_1,c_2,c_3)$.

\item[\bf ${\bold B_0}$.] $J_{B0}=J_{B0}(x_2,x_3):=\left[
           \begin{array}{ccc}
              x_2  & 1 & 0 \\
             0 & x_2 & 0 \\
             0 & 0 & x_3 \\
           \end{array}
         \right]$,
$F_{B0}=F_{B0}(x_2,x_3):=\left[
           \begin{array}{ccc}
             0 & 0 & x_2^2x_3 \\
             1 & 0 & -x_2^2-2x_2x_3 \\
             0 & 1 & 2x_2+x_3 \\
           \end{array}
         \right]$,\\
$P_{B0}=P_{B0}(x_2,x_3)=\left[\begin{array}{ccc}
{\frac{x_2x_3}{x_2-x_3}} & -{\frac{x_3(2\,x_2-x_3)}{(x_2-x_3)^2}} & {\frac{{x_2}^2}{(x_2-x_3)^2}}\\
\noalign{\medskip} -{\frac{x_2+x_3}{x_2-x_3}} & {\frac {2x_2}{(x_2-x_3)^2}} & -{\frac{2x_2}{(x_2-x_3)^2}}\\
\noalign{\medskip} \frac1{x_2-x_3} & -\frac1{(x_2-x_3)^2} & \frac1{(x_2-x_3)^2}
 \end{array}\right]$.
\item[\bf ${\bold B_1}$.] $J_{B1}:=J_{B0}(x_2,c_3)$, $F_{B1}:=F_{B0}(x_2,c_3)$, $P_{B1}=P_{B0}(x_2,c_3)$.
\item[\bf ${\bold B_2}$.] $J_{B2}:=J_{B0}(c_2,x_3)$, $F_{B2}:=F_{B0}(c_2,x_3)$, $P_{B2}=P_{B0}(c_2,x_3)$.
\item[\bf ${\bold B_3}$.] $J_{B3}:=J_{B0}(c_2,c_3)$, $F_{B3}:=F_{B0}(c_2,c_3)$, $P_{B3}=P_{B0}(c_2,c_3)$.

\item[\bf ${\bold C_0}$.] $J_{C0}=J_{C0}(x_2,x_3):=\left[
           \begin{array}{ccc}
             x_3  & 0 & 1 \\
             1 & x_3 & -x_2 \\
             0 & 0 & x_3 \\
           \end{array}
         \right]$,
$F_{C0}=F_{C0}(x_3):=\left[
           \begin{array}{ccc}
             0 & 0 & x_3^3 \\
             1 & 0 & -3x_3^2 \\
             0 & 1 & 3x_3       \\
           \end{array}
         \right]$,\\
$P_{C0}=P_{C0}(x_2,x_3)=\left[\begin{array}{ccc}
-x_3 & x_3^2 & x_2x_3+1\\
\noalign{\medskip} 1 & -2\,x_3 & -x_2\\
\noalign{\medskip} 0 & 1 & 0
 \end{array}\right]$.
\item[\bf ${\bold C_1}$.] $J_{C1}:=J_{C0}(1,c_3)$, $F_{C1}:=F_{C0}(c_3)$, $P_{C1}=P_{C0}(1,c_3)$.

\item[\bf ${\bold D_0}$.] $J_{D0}=J_{D0}(x_1,x_2,x_3):=\left[
           \begin{array}{ccc}
             x_1  & -x_2 & 0 \\
             x_2 & x_1 & 0 \\
             0 & 0 & x_3 \\
           \end{array}
         \right]$,\\
$F_{D0}=F_{D0}(x_1,x_2,x_3):=\left[
           \begin{array}{ccc}
             0 & 0 & (x_1^2+x_2^2)x_3 \\
             1 & 0 & -x_1^2-2x_1x_3-x_2^2 \\
             0 & 1 & 2x_1+x_3 \\
           \end{array}
         \right]$,\\
$P_{D0}=P_{D0}(x_1,x_2,x_3)=\left[\begin{array}{ccc}
-\frac{x_3(2\,x_1-x_3)}{(x_1-x_3)^2+x_2^2} & {\frac{x_3(x_1^2-x_1x_3-x_2^2)}{x_2((x_1-x_3)^2+x_2^2)}} &
{\frac{x_1^2+x_2^2}{(x_1-x_3)^2+x_2^2}}\\
\noalign{\medskip} {\frac{2x_1}{(x_1-x_3)^2+x_2^2}} & -{\frac{x_1^2-x_2^2-x_3^2}{x_2( (x_1-x_3)^2+x_2^2)}} &
-\frac{2x_1}{(x_1-x_3)^2+x_2^2}\\
\noalign{\medskip} -\frac1{(x_1-x_3)^2+x_2^2} & {\frac{x_1-x_3}{x_2((x_1-x_3)^2+x_2^2)}} & \frac1{(x_1-x_3)^2+x_2^2}
 \end{array}\right]$.
\item[\bf ${\bold D_1}$.] $J_{D1}:=J_{D0}(x_1,x_2,c_3)$, $F_{D1}:=F_{D0}(x_1,x_2,c_3)$, $P_{D1}=P_{D0}(x_1,x_2,c_3)$.
\item[\bf ${\bold D_2}$.] $J_{D2}:=J_{D0}(c_1,c_2,x_3)$, $F_{D2}:=F_{D0}(c_1,c_2,x_3)$, $P_{D2}=P_{D0}(c_1,c_2,x_3), c_2\not=0$.
\item[\bf ${\bold D_3}$.] $J_{D3}:=J_{D0}(c_1,c_2,c_3)$, $F_{D3}:=F_{D0}(c_1,c_2,c_3)$, $P_{D3}=P_{D0}(c_1,c_2,c_3), c_2\not=0$.
\end{itemize}

\end{document}